\documentclass{gtmon_a}
\pdfoutput=1

\proceedingstitle{Groups, homotopy and configuration spaces (Tokyo
2005)}
\conferencestart{5 July 2005}
\conferenceend{11 July 2005}
\conferencename{Groups, homotopy and configuration spaces, 
in honour of Fred Cohen's 60th birthday}
\conferencelocation{University of Tokyo, Japan}

\editor{Norio Iwase}
\givenname{Norio}
\surname{Iwase}

\editor{Toshitake Kohno}
\givenname{Toshitake}
\surname{Kohno}

\editor{Ran Levi}
\givenname{Ran}
\surname{Levi}

\editor{Dai Tamaki}
\givenname{Dai}
\surname{Tamaki}

\editor{Jie Wu}
\givenname{Jie}
\surname{Wu}

\title[Odd-primary exponents of Lie groups]{Odd-primary homotopy exponents of\\compact simple Lie groups}

\author{Donald M Davis}  
\givenname{Donald M}
\surname{Davis}
\address{Department of Mathematics\\Lehigh University\\\newline
Bethlehem, PA 18015\\USA}
\email{dmd1@lehigh.edu}

\author{Stephen D Theriault}
\givenname{Stephen D}
\surname{Theriault}
\address{Department of Mathematical Sciences\\University of Aberdeen\\\newline
Aberdeen AB24 3UE\\UK}
\email{s.theriault@maths.abdn.ac.uk}


\volumenumber{13}
\issuenumber{}
\publicationyear{2008}
\papernumber{08}
\startpage{195}
\endpage{201}
\doi{}
\MR{}
\Zbl{}
\received{25 October 2005}
\revised{14 September 2006}
\accepted{11 October 2006}
\published{22 February 2008}
\publishedonline{22 February 2008}
\proposed{}
\seconded{}
\corresponding{}

\arxivreference{math.AT/0601441}

\makeatletter
\def\cnewtheorem#1[#2]#3{\newtheorem{#1}{#3}[section]
\expandafter\let\csname c@#1\endcsname\c@subsection}

\AtBeginDocument{\let\bar\wbar\let\tilde\wtilde\let\hat\what}
\let\Bbb\mathbb
\makeop{SU}
\makeop{Spin}
\makeop{Sp}

\swapnumbers
\cnewtheorem{thm}[subsection]{Theorem}
\cnewtheorem{cor}[subsection]{Corollary}
\cnewtheorem{lem}[subsection]{Lemma}
\cnewtheorem{sublem}[subsection]{Sublemma}
\cnewtheorem{prop}[subsection]{Proposition}
\theoremstyle{remark}
\cnewtheorem{defin}[subsection]{Definition}
  \let\c@equation\c@subsection
  \def\theequation{\bf \thesubsection}
\makeatother

\def\mapright#1{\ \smash{\mathop{\longrightarrow}\limits^{#1}}\ }

\def\vp{v_1^{-1}\pi}

\def\on{\operatorname}

\def\a{\alpha}
\def\bz{{\Bbb Z}}

\def\tfrac{\textstyle\frac}

\def\equ{\begin{equation}}

\def\endeq{\end{equation}}

\def\coker{\on{coker}}

\def\exp{\on{exp}}

\def\lfl{\lfloor}
\def\rfl{\rfloor}

\def\({\bigg(}
\def\[{\bigg[}
\def\){\bigg)}
\def\]{\bigg]}
\def\colon{{:}\;}

\begin{document}

\begin{asciiabstract}
We note that a recent result of the second author yields upper bounds
for odd-primary homotopy exponents of compact simple Lie groups which
are often quite close to the lower bounds obtained from v_1-periodic
homotopy theory.

\end{asciiabstract}

\begin{htmlabstract}
We note that a recent result of the second author yields upper bounds
for odd-primary homotopy exponents of compact simple Lie groups which
are often quite close to the lower bounds obtained from v<sub>1</sub>&ndash;periodic
homotopy theory.
\end{htmlabstract}

\begin{abstract} 
We note that a recent result of the second author yields upper bounds
for odd-primary homotopy exponents of compact simple Lie groups which
are often quite close to the lower bounds obtained from $v_1$--periodic
homotopy theory.
\end{abstract}

\maketitle

\section{Statement of results}\label{intro}
The homotopy $p$--exponent of a topological space $X$, denoted $\exp_p(X)$, is the largest $e$
such that some homotopy group  $\pi_i(X)$ contains a $\bz/p^e$--summand.\footnote{Some authors
(eg \cite{Th1}) say that $p^e$ is the  homotopy $p$--exponent.} In work dating
back to 1989, the first author and collaborators have obtained lower bounds for $\exp_p(X)$ for all compact simple Lie groups
$X$ and all primes $p$ by using $v_1$--periodic homotopy theory. Recently, the second author \cite{Th1} proved a general result,
stated here as \fullref{thlem}, which can yield upper bounds for homotopy exponents of spaces which map to a sphere.
In this paper, we show that these two bounds often lead to a quite narrow range of values for $\exp_p(X)$
when $p$ is odd and $X$ is a compact simple Lie group.

Our first new result, which will be proved in \fullref{sect2},
combines \fullref{thlem} with a classical result of Borel and Hirzebruch.
\begin{thm}\label{SUupper}  Let $p$ be odd.

{\rm(a)}\qua If $n<p^2+p$, then $\exp_p(\SU(n))\le n-1+\nu_p((n-1)!)$.

{\rm(b)}\qua If $n\ge p^2+1$, then $\exp_p(\SU(n))\le n+p-3+\binom{\lfl\frac{n-2}{p-1}\rfl-p+2}2$.
\end{thm}
Here and throughout, $\nu_p(-)$ denotes the exponent of $p$ in an integer, $p$ is an odd prime, and $\lfl x\rfl$ denotes the integer part of $x$ . All spaces are
localized at $p$. It is useful to note 
the elementary fact that
$$\nu_p(m!)=\lfl\tfrac mp\rfl+\lfl\tfrac m{p^2}\rfl+\cdots,$$
and the well-known fact that $\nu_p(m!)\le\lfl\frac{m-1}{p-1}\rfl$.

\fullref{SUupper}(a) compares nicely with the following known result.
\begin{thm}\label{SUlower}\

{\rm(a)\qua (Davis and Sun \cite[1.1]{DS1})}\qua For any prime $p$, $\exp_p(\SU(n))\ge n-1+\nu_p(\lfl\tfrac{n}p\rfl!)$.

{\rm(b)\qua (Davis and Yang \cite[1.8]{DY})}\qua If $p$ is odd, $1\le t<p$, and $tp-t+2\le n\le tp+1$, then $\exp_p(\SU(n))\ge n$.
\end{thm}

Thus we have the following corollary, which gives the only values of $n>p$ in which the precise value of $\exp_p(\SU(n))$ is known.
\begin{cor} If $p$ is an odd prime, and $n=p+1$ or $n=2p$, then $\exp_p(\SU(n))=n$.\end{cor}

When $n=p+1$, this  was known (although perhaps never published)
since, localized at $p$, we have
$\SU(p+1)\simeq B(3,2p+1)\times S^5\times\cdots \times
S^{2p-1}$, the exponent of which follows from \fullref{B3} together with the result of Cohen, Moore, and Neisendorfer \cite{CMN} that if $p$ is odd, then $\exp_p(S^{2n+1})=n$. Here and throughout,  $B(2n+1,2n+1+q)$ denotes an $S^{2n+1}$--bundle over $S^{2n+1+q}$ with attaching map
$\a_1$ a generator of $\pi_{2n+q}(S^{2n+1})$, and $q=2p-2$. Note also that the result of \cite{CMN} implies that
if $n\le p$, then $\exp_p(\SU(n))=\exp_p(S^3\times\cdots\times
S^{2n-1})=n-1$.
\begin{prop}\label{B3} If $p$ is odd, then $\exp_p(B(3,2p+1))=p+1$, while if $n>1$, then $n+p-1\le\exp_p(B(2n+1,2n+1+q))\le n+p$.
\end{prop}
\begin{proof} 
This just combines Bendersky, Davis and Mimura \cite[1.3]{BDMi} for
the lower bound and Theriault \cite[2.1]{Th1} for the upper
bound.\end{proof}

Upper and lower bounds for the $p$--exponents of $\Sp(n)$ and $\Spin(n)$
can be extracted from Theorems \ref{SUupper} and \ref{SUlower} using
long-known relationships of their $p$--localizations to that of
appropriate $\SU(m)$.  Indeed, Harris \cite{Har} showed that there
are $p$--local equivalences
\begin{eqnarray}\label{SUdecomp}\SU(2n)&\simeq& \Sp(n)\times (\SU(2n)/\Sp(n))\\
\Spin(2n+1)&\simeq& \Sp(n)\label{Har2}\\
\Spin(2n+2)&\simeq& \Spin(2n+1)\times S^{2n+1}.\label{Har3}\end{eqnarray}
Combining this with Theorems \ref{SUupper} and \ref{SUlower} leads to the following corollary.
\begin{cor} Let $p$ be odd.
\begin{enumerate}
\item $\exp_p(\Spin(2n+2))=\exp_p(\Spin(2n+1))=\exp_p(\Sp(n))\le\exp_p(\SU(2n))$, which is bounded according to
\fullref{SUupper}.

\item $\exp_p(\Sp(n))\ge 2n-1+\nu_p(\lfl\frac{2n}p\rfl!)$.

\item If $1\le t<p$, and $tp-t+2\le 2n\le tp+1$, then $\exp_p(\Sp(n))\ge 2n$.
\end{enumerate}
\end{cor}
\begin{proof} 
The second and third parts of (1) are immediate from \eqref{Har2} and
\eqref{SUdecomp}, while the first equality of (1) follows from
\eqref{Har3} and the fact that $\exp_p(\Spin(2n+1))\ge
\exp_p(S^{2n+1})$, which is a consequence of part (2) and
\eqref{Har2}.  For parts (2) and (3), we need to know that the
homotopy classes yielding the lower bounds for $\exp_p(\SU(2n))$ given
in \fullref{SUlower} come from its $\Sp(n)$ factor in \eqref{SUdecomp}.
To see this, we first note that in Bendersky and Davis
\cite[1.2]{BDspin} it was proved that, if $p$ is odd and $k$ is odd,
then
\begin{equation}\label{SUSp}\vp_{2k}(\Sp(n);p)\approx\vp_{2k}(\SU(2n);p).
\end{equation}
These denote the $p$--primary $v_1$--periodic homotopy groups, which appear as summands of actual homotopy groups.
The proofs of \cite[1.1]{DS1} and \cite[1.8]{DY}, which yielded \fullref{SUlower}, were obtained by computing
 $\vp_{2k}(\SU(n);p)$ for certain $k\equiv n-1$ mod 2. When applied to $\SU(2n)$, these groups are in $\vp_{2k}(\SU(2n);p)$ with $k$ odd, and so by \eqref{SUSp} they appear in the $\Sp(n)$ factor. \end{proof}
 
For all $(X,p)$ with $X$ an exceptional Lie group and $p$ an odd
prime, except $(E_7,3)$ and $(E_8,3)$, we can make an excellent
comparison of bounds for $\exp_p(X)$ using results in the literature.
We use splittings of the torsion-free cases tabulated in
\cite[1.1]{BDMi}, but known much earlier (Mimura, Nishida and Toda
\cite{MNT}).  In \fullref{exctabl}, we list the range of possible
values of $\exp_p(X)$ when the precise value is not known.  We also
list the factor in the product decomposition which accounts for the
exponent. Finally, in cases in which the exponent bounds do not follow
from results already discussed, we provide references.  Here
$B(n_1,\ldots,n_r)$ denotes a space built from fibrations involving
$p$--local spheres of the indicated dimensions and equivalent to a
factor in a $p$--localizaton of a special unitary group or quotient of
same. Also, $B_2(3,11)$ denotes a sphere-bundle with attaching map
$\a_2$, and $W$ denotes a space constructed by Wilkerson and shown in
\cite[1.1]{Th2} to fit into a fibration $\Omega K_5\to B(27,35)\to
W$. Finally, $K_3$ and $K_5$ denotes Harper's space as described in
Bendersky and Davis \cite{BDF4} and Theriault \cite{Th1}.

\begin{thm}\label{excthm} The homotopy $p$--exponents of exceptional
  Lie groups are as given in \mbox{\fullref{exctabl}}.\end{thm}

\begin{table}[ht!]
\begin{center}
\renewcommand{\arraystretch}{1.15}
\caption{Homotopy exponents of exceptional Lie groups}\label{exctabl}
\begin{tabular}{cc|ccc}
$X$&$p$&$\exp_p(X)$&Factor&Reference\\
\hline
$G_2$&$3$&$6$&$B_2(3,11)$&\cite[1.3]{BDMi}, \cite[2.2]{Th1}\\
$G_2$&$5$&$6$&$B(3,11)$&\\
$G_2$&$>5$&$5$&$S^{11}$&\\
\hline
$F_4,E_6$&$3$&$12$&$K_3$&\cite[1.6]{BDF4}, \cite[1.2]{Th1}\\
$F_4,E_6$&$5,7$&$11,12$&$B(23-q,23)$&\\
$F_4,E_6$&$11$&$12$&$B(3,23)$&\\
$F_4,E_6$&$>11$&$11$&$S^{23}$&\\
\hline
$E_7$&$5$&$18,19,20$&$B(3,11,19,27,35)$&factor of $\SU(18)$\\
$E_7$&$7$&$17,18,19$&$B(11,23,35)$&factor of $\SU(18)$\\
$E_7$&$11,13$&$17,18$&$B(35-q,35)$&\\
$E_7$&$17$&$18$&$B(3,35)$&\\
$E_7$&$>17$&$17$&$S^{35}$&\\
\hline
$E_8$&$5$&$30,31$&$W$&\cite[1.1]{Rep}, \cite[1.2]{Th2}\\
$E_8$&$7$&$29,30,31,32$&$B(23,35,47,59)$&\cite[1.4]{BDMi}, \fullref{onecase}\\
$E_8$&$11-23$&$29,30$&$B(59-q,59)$&\\
$E_8$&$29$&$30$&$B(3,59)$&\\
$E_8$&$>29$&$29$&$S^{59}$&
\end{tabular}
\end{center}
\end{table}

\section[Proof of \ref{SUupper}]{Proof of \fullref{SUupper}}\label{sect2}
In \cite[Lemma 2.2]{Th1}, the second author proved the following result.
\begin{lem}{\rm\cite[2.2,2.3]{Th1}}\qua Suppose there is a homotopy fibration
$$F\to E\mapright{q} S^{2n+1}$$
where $E$ is simply-connected or an $H$--space and $$|\coker(\pi_{2n+1}(E)\mapright{q_*}\pi_{2n+1}(S^{2n+1}))|\le p^r.$$
Then $\exp_p(E)\le r+\max(\exp_p(F),n)$.\label{thlem}\end{lem}
In \cite[2.2]{Th1}, it was required that $E$ be an $H$--space, but  \cite[2.3]{Th1}  noted that
if $E$ is not an $H$--space, the desired conclusion can be obtained by applying the loop-space functor to
the fibration. We require $E$ to be simply-connected so that we do not loop away a large fundamental group. 
We now use this lemma to prove \fullref{SUupper}.

\begin{proof}[Proof of \fullref{SUupper}] The proof is by induction on $n$.
Let the odd prime $p$ be implicit, and let $\SU'(n)$ denote the factor in the $p$--local product decomposition \cite{MNT}
of $\SU(n)$ which is built from spheres of dimension congruent to $2n-1$ mod $2p-2$.  By the induction hypothesis, the exponents of the other factors are $\le$
the asserted amount.
We will apply \fullref{thlem} to the fibration
$$\SU'(n-p+1)\to \SU'(n)\mapright{q} S^{2n-1}.$$
In order to determine $|\coker(\pi_{2n-1}(\SU'(n))\mapright{q_*}\pi_{2n-1}(S^{2n-1}))|$, we use
the classical result of Borel and Hirzebruch (\cite[26.7]{BH}) that $$\pi_{2n-2}(\SU(n-1))\approx\bz/(n-1)!.$$
When localized at $p$, it is clear that its $p$--component $\bz/p^{\nu_p((n-1)!)}$ must come from the $\SU'(n-p+1)$--factor
in the product decomposition of $\SU(n-1)$,
since\break $\pi_{2n-2}(\SU(n-1))$ is built from the classes $\a_{i}\in\pi_{2n-2}(S^{2n-1-iq})_{(p)}$.
Thus
$$\pi_{2n-2}(\SU'(n-p+1))\approx\bz/p^{\nu_p((n-1)!)},$$
and the exact sequence
$$\pi_{2n-1}(\SU'(n))\mapright{q_*}\pi_{2n-1}(S^{2n-1})\to\pi_{2n-2}(\SU'(n-p+1))$$
implies \begin{equation}\label{cok}\nu_p(|\coker(q_*)|)\le\nu_p((n-1)!).\end{equation}

(a)\qua By the induction hypothesis, $\exp_p(\SU'(n-p+1))\le n-p+\nu_p((n-p)!)$. By hypothesis, $n-p<p^2$ and hence $\nu_p((n-p)!)\le p-1$.
Thus $\exp_p(\SU'(n-p+1))\le n-1$, and so by \fullref{thlem} and \eqref{cok}
$$\exp_p(\SU'(n))\le \nu_p(|\coker(q_*)|)+ n-1\le\nu_p((n-1)!)+n-1,$$
as claimed.

(b)\qua By (a), part (b) is true if $p^2+1\le n\le p^2+p-1$.
Let $n\ge p^2+p$, and assume the theorem is true for $\SU'(n-p+1)$. 
Then by \fullref{thlem} and the induction hypothesis
$$\exp_p(\SU'(n))\le\nu((n-1)!)+n-p+1+p-3+\binom{\lfl\frac{n-p-1}{p-1}\rfl-p+2}2.$$
Note that even if $\exp_p(\SU'(n-p+1))$ happened to be less than $n-1$, our upper bound for it is $\ge n-1$,
and so this bound for $\exp_p(\SU'(n))$ is still a correct deduction from \fullref{thlem}.

Since $\nu_p((n-1)!)\le\lfl\frac{n-2}{p-1}\rfl$, we obtain
\begin{eqnarray*}\exp_p(\SU'(n))
&\le&\left\lfl\frac{n-2}{p-1}\right\rfl+n-2+\binom{\lfl\frac{n-2}{p-1}\rfl-p+1}2\\
&=&\left\lfl\frac{n-2}{p-1}\right\rfl+n-2+\binom{\lfl\frac{n-2}{p-1}\rfl-p+2}2-\binom{\lfl\frac{n-2}{p-1}\rfl-p+1}1\\
&=&n+p-3+\binom{\lfl\frac{n-2}{p-1}\rfl-p+2}2,\end{eqnarray*}
as desired.
\end{proof}

The result in part (b) could be improved somewhat by a more delicate numerical argument.

Part (b) of the following result was used in \fullref{exctabl}.
\begin{prop}\label{onecase} Let $p=7$.

{\rm(a)}\qua $\exp_7(B(23,35,47))\le 25$.

{\rm(b)}\qua $\exp_7(B(23,35,47,59))\le 32$.\end{prop}
\begin{proof} The thing that makes this require special attention is that these spaces are not a factor
of an $\SU(n)$, because they do not contain an $S^{11}$. There are fibrations
$$B(23,35)\to  B(23,35,47)\to  S^{47}$$
$$B(23,35,47)\to B(23,35,47,59)\to S^{59}.\leqno{\hbox{and}}$$
Since, localized at 7, $\pi_{46}(S^{23})\approx\pi_{46}(S^{35})\approx\bz/7$, we have $|\pi_{46}(B(23,35))|\le7^2$,
and similarly $|\pi_{58}(B(23,35,47))|\le7^3$. (In fact, it is easily seen that these are cyclic groups of the
indicated order.) Using \fullref{thlem} and that $\exp_7(B(23,35))\le18$ by \fullref{B3}, we obtain 
$$\exp_7(B(23,35,47))\le 2+\max(18,23)=25,$$
and then
$$\exp_7(B(23,35,47,59))\le 3+\max(25,29)=32.\proved$$
\end{proof}

\bibliographystyle{gtart}
\bibliography{link}

\end{document}